\documentclass{amsart}
\usepackage{amsmath}
\usepackage{amsthm}
\usepackage{amssymb}
\usepackage{epsfig}

\theoremstyle{plain}
\newtheorem{thm}{Theorem}[section]
\newtheorem{lem}[thm]{Lemma}
\newtheorem{prop}[thm]{Proposition}

\theoremstyle{definition}
\newtheorem{example}[thm]{Example}
\newtheorem{defn}[thm]{Definition}
\theoremstyle{remark}

\newcommand{\bm}[1]{\mbox{\boldmath $#1$}}
\newcommand{\ul}[1]{\underline{#1}}
\newcommand{\B}{\mathfrak{S}_n^B}
\newcommand{\w}{w_{0}^{B}}
\newcommand{\X}{X_n^B}
\newcommand{\Y}{Y_n^B}

\title{On Expected Factors in Reduced Decompositions in Type $B$}
\author{Bridget Eileen Tenner}
\address{Department of Mathematics, Massachusetts Institute of Technology, 77 Massachusetts Ave, Cambridge, MA 02139, USA}
\email{bridget@math.mit.edu}
\date{August 10, 2005}

\begin{document}

\begin{abstract}

The expected number of Yang-Baxter moves appearing in a reduced decomposition of the longest element of the Coxeter group of type $B_n$ is computed to be $2-4/n$.  For the same element, the expected number of $0101$ or $1010$ factors appearing in a reduced decomposition is $2/(n^2-2)$.

\end{abstract}

\maketitle

\section{Introduction}

Let $\B$ denote the Coxeter group of type $B$ on $\{\pm 1, \ldots, \pm n\}$, also known as the hyperoctahedral group.  This is the group of signed permutations $w$ on $\{\pm 1, \ldots, \pm n\}$, with the requirement that $w(-i) = -w(i)$.  For ease of notation, a negative sign will be written beneath an integer, rather than before it.  An element $w \in \B$ will be written $w = w(1) \cdots w(n)$, where $w$ maps $i \mapsto w(i)$ (and $\ul{i} \mapsto \ul{w(i)}$).  Notice that the element $w$ is completely determined by the values of $w(1), \ldots, w(n)$.  The group $\B$ is generated by the transpositions $\{s_i : 0 \le i \le n-1\}$, which satisfy the Coxeter relations:

\begin{eqnarray}
& s_i^2 = 1 &\text{ for all }i; \label{s_i}\\
& s_is_j = s_j s_i &\text{ if } |i-j| > 1; \label{s_is_j}\\
& s_is_{i+1}s_i = s_{i+1}s_is_{i+1} & \text{ for } 1 \le i \le n-2; \label{longbraid}\\
& s_0s_1s_0s_1 = s_1s_0s_1s_0. \label{01}&
\end{eqnarray}

A map is written to the left of its input, meaning that for $i > 0$, the product $s_iw$ changes the positions of the values $i$ and $i+1$ (and $\ul{i}$ and $\ul{i+1}$) in $w$, whereas $ws_i$ changes the values in positions $i$ and $i+1$ in $w$.  The transposition $s_0$ changes the sign of the first entry, so $s_0w$ interchanges $\ul{1}$ and $1$ in $w$, and $ws_0 = \ul{w(1)}w(2)\cdots w(n)$.  An extensive treatment of Coxeter groups can be found in \cite{bjornerbrenti}.

Every element in $\B$ can be written as a product of the transpositions $\{s_i: 0 \le i \le n-1\}$.  The minimum number of transpositions required for a product to equal $w$ is the \emph{length} of $w$, denoted $\ell(w)$.  The \emph{longest element} in $\B$ is $\w = \ul{1} \ul{2} \cdots \ul{n}$, and $\ell(\w) = n^2$.

\begin{defn}
A \emph{reduced decomposition} for an element $w$ of length $\ell = \ell(w)$ is a string $i_1 \cdots i_{\ell}$ such that $w = s_{i_1} \cdots s_{i_{\ell}}$.  The set $R(w)$ consists of all reduced decompositions of $w$.
\end{defn}

\begin{defn}
A consecutive substring of a reduced decomposition is a \emph{factor}.  For any $j>0$, a factor of the form $j (j + 1) j$ or $(j+1)j(j+1)$ in a reduced decomposition will be called a \emph{Yang-Baxter move}, and a factor $0101$ or $1010$ will be called a \emph{01 move}.
\end{defn}

The symmetric group $\mathfrak{S}_n$ of unsigned permutations is the Coxeter group of type $A$.  It is generated by the transpositions $\{s_i: 1 \le i \le n-1\}$, which are subject to the relations in equations~\eqref{s_i}-\eqref{longbraid}.  The longest element in $\mathfrak{S}_n$ is $w_0 = n(n-1) \cdots 1$, which has length $\binom{n}{2}$.  In \cite{reiner}, Reiner computes the following somewhat surprising result.

\begin{thm}[Reiner]\label{reinerthm}
The expected number of Yang-Baxter moves in a reduced decomposition of $w_0 \in \mathfrak{S}_n$ is $1$ for all $n \ge 3$.
\end{thm}

This paper presents results for Coxeter groups of type $B$ that are analogous to Theorem~\ref{reinerthm}.   In type $A$, the expectation of factors corresponding to the Coxeter relation in equation~\eqref{longbraid} was computed.  For the hyperoctahedral group, factors corresponding to the Coxeter relations in each of equations~\eqref{longbraid} and~\eqref{01} will be treated.
Theorem~\ref{E(X)thm} calculates that the expected number of Yang-Baxter moves in a reduced decomposition of $\w \in \B$ is $2 - 4/n$, and Theorem~\ref{E(Y)thm} shows that the expected number of $01$ moves is $2/(n^2 - 2)$.  Unlike Reiner's result, both of these expectations are dependent upon $n$.  Moreover, in the context of Theorem~\ref{reinerthm}, the value $2 - 4/n$ seems quite plausible since the length of $\w \in \B$ is approximately twice that of $w_0 \in \mathfrak{S}_n$.

A variety of tools are used to prove Theorems~\ref{E(X)thm} and~\ref{E(Y)thm}, several of which are discussed in Section 2.  Section 3 computes the expected number of Yang-Baxter moves in elements of $R(\w)$, and Section 4 does the same for $01$ moves.

\section{Vexillary elements, shapes, and hook lengths in type $B$}

In \cite{stanley}, Stanley shows that for a vexillary element $v \in \mathfrak{S}_n$,
\begin{equation}\label{vexA}
\#R(v) = f^{\lambda(v)},
\end{equation}
\noindent where $f^{\lambda(v)}$ is the number of standard Young tableaux of a particular shape $\lambda(v)$.  This result is central to the proof of Theorem~\ref{reinerthm}.

There are numerous definitions of vexillarity in type $A$, for example see \cite{macdonald} and \cite{tenner}.  One definition is that $w \in \mathfrak{S}_n$ is vexillary if it is $2143$-avoiding.  Billey and Lam define a notion of vexillary for type $B$ in \cite{billeylam}.  Their definition is in terms of Stanley symmetric functions and Schur $Q$-functions, and they prove its equivalence with a statement about pattern avoidance, now of signed permutations.  This latter statement will be given as the definition here, and it follows from the work of Billey and Lam that it generalizes equation~\eqref{vexA} to type $B$ in the appropriate way.

\begin{defn}
An element $w \in \B$ is \emph{vexillary for type $B$} if $w = w(1) \cdots w(n)$ avoids the following patterns:
\begin{center}
$\begin{array}{lll}
21 & \ul{3}2\ul{1} & 2\ul{3}4\ul{1}\\
\ul{2}\ul{3}4\ul{1} & 3\ul{4}\ul{1}\ul{2} & \ul{3}\ul{4}1\ul{2}\\
\ul{3}\ul{4}\ul{1}\ul{2} & \ul{4}1\ul{2}3 & \ul{4}\ul{1}\ul{2}3
\end{array}$
\end{center}
\end{defn}

Note that patterns in signed permutations must maintain their signs.  For example, $\ul{1}\ul{2}$ is not an instance of the pattern $21$, even though $\ul{1} > \ul{2}$.

\begin{example}
$2\ul{1}\ul{4}3 \in \mathfrak{S}_4^B$ is vexillary for type $B$, but $2\ul{1}43 \in \mathfrak{S}_4^B$ is not.
\end{example}

To each element $w \in \B$, Billey and Lam define a shifted shape $\lambda^B(w)$ as follows.

\begin{defn}\label{shapealg}
Let $w = w(1) \cdots w(n) \in \B$.
\begin{enumerate}
\item Write $\{w(1), \ldots, w(n)\}$ in increasing order and call this $u \in \B$.
\item Let $v \in \mathfrak{S}_n$ be the (vexillary) permutation $u^{-1}w$.
\item Let $\mu$ be the partition with (distinct) parts $\{|u_i| : u_i < 0\}$.
\item Let $U$ be any standard shifted Young tableau of shape $\mu$, and let $V$ be any standard Young tableau whose shape is the transpose of the partition with parts $\{c_1, \ldots, c_n\}$, where
\begin{equation*}
c_i = \#\{j : j>i \text{ and } v(j) < v(i)\}.
\end{equation*}
\item Embed $U$ in the shifted shape $\delta = (n, n-1, \ldots, 1)$.
\item Fill in the rest of $\delta$ with $1', \ldots, k'$ starting from the rightmost column and labeling each column from bottom to top.  This gives the tableau $R$.
\item Obtain $S$ by adding $|\mu|$ to each entry of $V$, and glue $R$ to the left side of $S$ to obtain $T$.
\item Delete the box containing $1'$ from $T$.  If the remaining tableau is not shifted, apply jeu de taquin to fill in the box.  Do likewise for the box containing $2'$, then $3'$, and so on, stopping after procedure for the box containing $k'$.
\end{enumerate}
\noindent The (shifted) shape of the resulting tableau is $\lambda^B(w)$.
\end{defn}

\begin{example}
Suppose $w = 2\ul{1}\ul{4}3 \in \mathfrak{S}_4$.  Then $u = \ul{4}\ul{1}23$, $v = 3214$, $\mu = (4,1)$, and the tableau $V$ has shape $(2,1)$.  Five boxes of $\delta$ will be filled by primed numbers, and the final tableau has shifted shape $\lambda^B(w) = (6,2)$.
\end{example}

\begin{prop}[Billey-Lam]\label{redtab}
If $w \in \B$ is vexillary for type $B$, then
\begin{equation}\label{vexB}
\# R(w) = f^{\lambda^B(w)}
\end{equation}
\noindent where $f^{\lambda^B(w)}$ is the number of standard tableaux of shifted shape $\lambda^B(w)$.
\end{prop}

Equation~\eqref{vexB} will play an analogous role in the proofs of this paper to that played by equation~\eqref{vexA} in \cite{reiner}.  Hooks and hook-lengths for shifted shapes will also be important tools, as they facilitate the calculation of $f^{\lambda^B(w)}$.  Recall the hook-length formula for straight shapes (see \cite{ec2} for a more extensive treatment).

\begin{prop}
For a shape $\lambda \vdash N$,
\begin{equation*}
f^\lambda = \frac{N!}{\prod_{u \in \lambda} h(u)},
\end{equation*}
\noindent where $h(u)$ is the number of squares in $\lambda$ that are in the same column as $u$ but no higher, or in the same row as $u$ but no farther to the left.
\end{prop}

There is an analogous formula for shifted shapes (for more information, see \cite{sagan}).

\begin{prop}\label{Bhook}
For a shifted shape $\lambda^B \vdash N$,
\begin{equation*}
f^{\lambda^B} = \frac{N!}{\prod_{u \in \lambda^B} h^B(u)},
\end{equation*}
\noindent where $h^B(u)$ is the total number of the squares in $\lambda^B$ that are
\begin{enumerate}
\item In the same column as $u$ but no higher;
\item In the same row as $u$ but no farther to the left; or
\item In the $(k+1)$st row of $\lambda^B$ if $u$ is in the $k$th column of $\lambda^B$.
\end{enumerate}
\end{prop}

\begin{figure}[htbp]
\centering
\epsfig{file=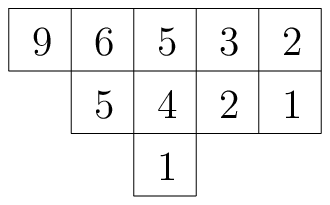, width=1in}
\caption{Hook lengths for $\lambda^B = (5,4,1)$.  $f^{\lambda^B} = 56$.}
\end{figure}

The final preliminary to proving the main results of this paper is the following lemma.  The proof is straightforward and omitted here.

\begin{lem}\label{rotate}
For $\w \in \B$ and $i \in [0,n-1]$,
\begin{equation*}
s_i\w s_i = \w.
\end{equation*}
\end{lem}

This indicates a $\mathbb{Z}/n^2\mathbb{Z}$-action on the set $R(\w)$ defined by
\begin{equation*}
s_{i_1}s_{i_2} \cdots s_{i_{n^2}} \mapsto s_{i_2} \cdots s_{i_{n^2}}s_{i_1}.
\end{equation*}

As with the other machinery discussed in this section, Lemma~\ref{rotate} has an analogous (though not identical) statement in type $A$ which is used in \cite{reiner}.

\section{Expectation of Yang-Baxter moves}

Consider the set $R(\w)$ with uniform probability distribution.  Let $\X$ be the random variable on reduced decompositions of $\w \in \B$ which counts the number of Yang-Baxter moves.

\begin{thm}\label{E(X)thm}
For all $n \ge 3$, $E(\X) = 2 - 4/n$.
\end{thm}

\begin{proof}

Fix $n \ge 3$.  For $k > 0$, let $\X[j,k]$ be the indicator random variable which determines whether the factor $i_j i_{j+1} i_{j+2}$ in a reduced decomposition $i_1 \cdots i_{n^2}$ in $R(\w)$ is of either form $k(k+1)k$ or $(k+1)k(k+1)$.  Therefore
\begin{equation*}
E(\X) = \sum_{j=1}^{n^2-2} \ \sum_{k=1}^{n-2} E(\X[j,k]).
\end{equation*}

\noindent The variables $\X[j,k]$ and $\X[j',k]$ have the same distribution by Lemma~\ref{rotate}, so in fact
\begin{equation*}
E(\X) = (n^2 - 2)\sum_{k=1}^{n-2} E(\X[1,k]).
\end{equation*}

If $\X[1,k](\bm{i}) = 1$ for $\bm{i} = i_1 \cdots i_{n^2} \in R(\w)$, then
\begin{equation*}
\bm{i} \in \{k(k+1)ki_4 \cdots i_{n^2},(k+1)k(k+1)i_4 \cdots i_{n^2}\}.
\end{equation*}
\noindent In both cases, $i_4 \cdots i_{n^2}$ is a reduced decomposition of
\begin{equation*}
w_k := s_ks_{k+1}s_k\w = \ul{1}\cdots \ul{(k-1)}\ul{(k+2)}\ul{(k+1)}\ul{(k)}\ul{(k+3)} \cdots \ul{n}.
\end{equation*}

\noindent Notice that $w_k$ is vexillary for type $B$ for all $k$.  Therefore, by Proposition~\ref{redtab},
\begin{equation}\label{E(X)3}
E(\X) = 2(n^2 - 2) \sum_{k=1}^{n-2} \frac{\#R(w_k)}{\#R(\w)} = 2(n^2 - 2) \sum_{k=1}^{n-2} \frac{f^{\lambda^B(w_k)}}{f^{\lambda^B(\w)}}
\end{equation}

The shifted shapes $\lambda^B(\w)$ and $\lambda^B(w_k)$ are easy to determine, as the signed permutation $u$ in Definition~\ref{shapealg} is $\ul{n} \cdots \ul{1}$ in both cases, so no boxes contain primed entries in the shifted tableau $T$.  Thus, the shifted shapes are
\begin{eqnarray}
\lambda^B(\w) &=& (2n-1,2n-3, \ldots, 3,1) \text{ and}\label{lambda(w)}\\
\lambda^B(w_k) &=& (2n-1,2n-3, \ldots, 2k+5,2k+1,2k,2k-1, \ldots, 3,1)\nonumber.
\end{eqnarray}

Recall the hook-length formula of Proposition~\ref{Bhook}, particularly the definition of the hooks $h^B$ in shifted shapes.  The only hook-lengths that do not cancel in the ratio $f^{\lambda^B(w_k)}/f^{\lambda^B(\w)}$ are as indicated in Figure~\ref{shapes}.

\begin{figure}[htbp]
\begin{center}
$\begin{array}{c@{\hspace{.2in}}c}
\epsfig{file=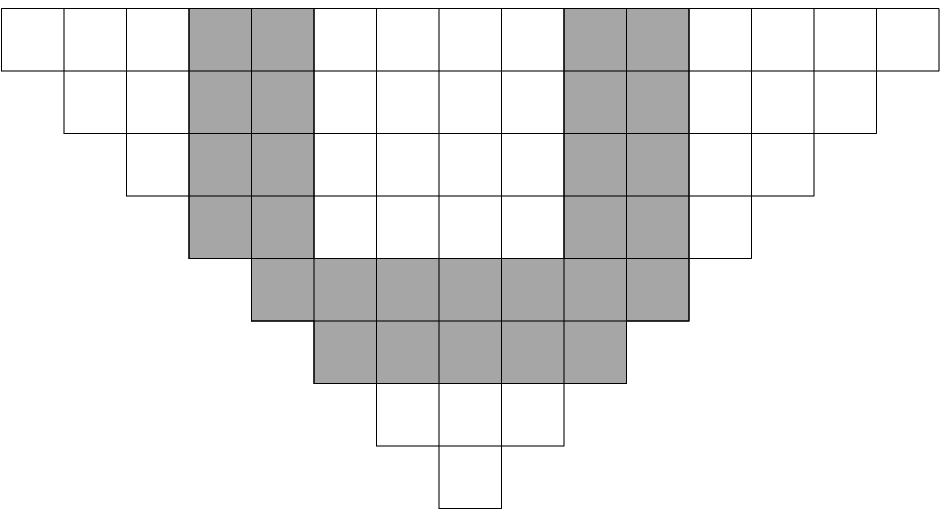, width=2in} & \hspace{.2in} \epsfig{file=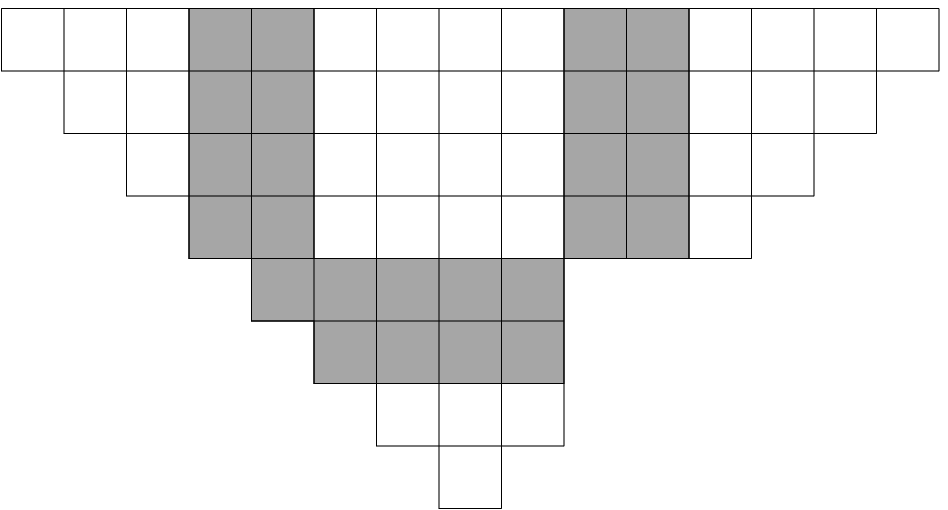, width=2in}
\end{array}$
\end{center}
\caption{The shifted shapes $\lambda^B(\w)$ and $\lambda^B(w_k)$ for $n = 8$ and $k = 2$.  The shaded boxes are where the hook-lengths are unequal.}\label{shapes}
\end{figure}

Consequently, equation~\eqref{E(X)3} can be written as
\begin{equation}\label{E(X)4}
E(\X) = \frac{1}{3}\binom{n^2}{2}^{-1}\ \sum_{k=1}^{n-2} C_k,
\end{equation}
\noindent where
\begin{multline*}
C_k = \frac{3\cdot 5\cdots (2k+3)}{2\cdot 4 \cdots (2k)} \cdot \frac{3\cdot 5 \cdots (2n-2k-1)}{2\cdot 4 \cdots (2n-2k-4)} \cdot \frac{(2k+4)(2k+6)\cdots (4k+4)}{(2k+1)(2k+3)\cdots(4k+1)}\\
 \cdot \frac{(4k+8)(4k+10)\cdots(2n+2k+2)}{(4k+5)(4k+7)\cdots(2n+2k-1)},
\end{multline*}
\noindent and empty products are defined to be $1$.

Notice that
\begin{equation*}
\frac{C_{k+1}}{C_k} = \frac{(2k+3)(4k+7)(2k+1)(n-k-2)(n+k+2)}{(4k+3)(k+2)(2n+2k+1)(2n-2k-1)(k+1)}
\end{equation*}
\noindent is a rational function in $k$.  Therefore, $\sum_{k=1}^{n-2} C_k$ is a hypergeometric series.  Following the notation in \cite{roy}, equation~\eqref{E(X)4} can be rewritten as
\begin{equation*}
E(\X) = \frac{1}{3}\binom{n^2}{2}^{-1}C_0\left(_5F_4\begin{pmatrix} 3/2, 7/4, 1/2, 2-n, 2+n\\ 3/4, 2,1/2+n, 1/2-n\end{pmatrix} - 1 \right).
\end{equation*}

The hypergeometric series in question can be computed via Dougall's theorem, as discussed in \cite{slater}.  The theorem states that
\begin{multline*}
_5F_4\begin{pmatrix} a, 1 + a/2, b, c, d\\a/2, 1+a-b, 1+a-c, 1+a-d\end{pmatrix}\\ = \frac{\Gamma(1+a-b)\Gamma(1+a-c)\Gamma(1+a-d)\Gamma(1+a-b-c-d)}{\Gamma(1+a)\Gamma(1+a-b-c)\Gamma(1+a-b-d)\Gamma(1+a-c-d)}.
\end{multline*}

\noindent It is not immediately obvious that Dougall's theorem applies to this particular series because of a potential pole.  However, the theorem does show that
\begin{eqnarray*}
_5F_4\begin{pmatrix} 3/2, 7/4, 1/2, 2-n, 2+x\\ 3/4, 2,1/2+n, 1/2-x\end{pmatrix} &=& \frac{\Gamma(2)\Gamma(1/2+n)\Gamma(1/2-x)\Gamma(n-x-2)}{\Gamma(5/2)\Gamma(n)\Gamma(-x)\Gamma(n-x-3/2)}\\
&=& \frac{(-x)_{n-2} (5/2)_{n-2}}{(n-1)!(1/2-x)_{n-2}}.
\end{eqnarray*}
\noindent Therefore there is no pole in this situation.  Letting $x$ approach $n$ shows that the desired hypergeometric series has sum $n/2$.

Finally,
\begin{equation*}
C_0 = 3 \cdot \frac{3 \cdot 5 \cdots (2n-1)}{2 \cdot 4 \cdots (2n-4)} \cdot \frac{4}{1} \cdot \frac{8 \cdot 10 \cdots (2n+2)}{5 \cdot 7 \cdots (2n-1)} = 6n(n^2-1),
\end{equation*}
\noindent which completes the proof:
\begin{equation*}
E(\X) = \frac{1}{3}\binom{n^2}{2}^{-1}6n(n^2-1)(n/2 - 1) = 2 - 4/n.
\end{equation*}

\end{proof}

As suggested earlier, it is appropriate that Yang-Baxter moves are approximately twice as common in elements of $R(\w)$ as in elements of $R(w_0)$, as 
\begin{equation*}
\ell(\w) = n^2 \approx 2\binom{n}{2} = 2\ell(w_0).
\end{equation*}

\section{Expectation of $01$ moves}

As in the previous section, consider the set $R(\w)$ with uniform probability distribution.  Let $\Y$ be the random variable on reduced decompositions of $\w \in \B$ which counts the number of $01$ moves.

\begin{thm}\label{E(Y)thm}
For all $n \ge 2$, $E(\Y) = 2/(n^2 - 2)$.
\end{thm}

\begin{proof}

As in the proof of Theorem~\ref{E(X)thm}, let $\Y[j]$ be the indicator random variable which determines whether the factor $i_ji_{j+1}i_{j+2}i_{j+3}$ in a reduced decomposition $i_1 \cdots i_{n^2}$ in $R(\w)$ is of either form $0101$ or $1010$.  Lemma~\ref{rotate} similarly applies, so
\begin{equation*}
E(\Y) = \sum_{j=1}^{n^2-3}E(\Y[j]) = (n^2 - 3)E(\Y[1]).
\end{equation*}

If $\Y[1](\bm{i}) = 1$ for $\bm{i} = i_1 \cdots i_{n^2} \in R(\w)$, then $\bm{i}$ is either $0101i_5 \cdots i_{n^2}$ or $1010i_5 \cdots i_{n^2}$.  The string $i_5 \cdots i_{n^2}$ is a reduced decomposition of
\begin{equation*}
w' = 12\ul{3}\ul{4}\cdots\ul{n}
\end{equation*}
\noindent in both situations.  As with $w_k$, the signed permutation $w'$ is vexillary for type $B$.  Therefore
\begin{equation*}
E(\Y) = 2(n^2 - 3) \frac{\#R(w')}{\#R(\w)} = 2(n^2 - 3) \frac{f^{\lambda^B(w')}}{f^{\lambda^B(\w)}}.
\end{equation*}

The shifted shape $\lambda^B(\w)$ is as in equation~\eqref{lambda(w)}.  Applying Definition~\ref{shapealg} to $w'$ proceeds as follows:

\begin{enumerate}
\item The signed permutation $u$ is $\ul{n} \cdots \ul{3}12$.
\item The vexillary permutation $v$ is $(n-1)n(n-2)(n-3) \cdots 21$.
\item The partition $\mu$ is $(n,n-1,\ldots, 4,3)$.
\item The shifted tableau $U$ has shape $(n,n-1, \ldots, 4,3)$ and the straight tableau $V$ has shape $(n-1,n-2,\ldots,3,2)$.
\end{enumerate}

\begin{figure}[htbp]
\centering
\epsfig{file=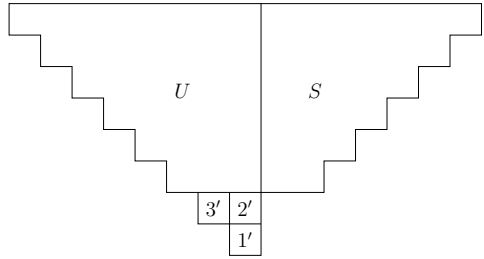, width=2.5in}
\caption{Step (7) of Definition~\ref{shapealg} applied to $w'$.}\label{shape(w')}
\end{figure}

Unlike in the cases of $\w$ or $w^k$ in the proof of Theorem~\ref{E(X)thm}, there will be boxes of $T$ containing primed numbers, specifically $1'$, $2'$, and $3'$, as in Figure~\ref{shape(w')}.  However, removing $1'$ leaves a shifted tableau so jeu de taquin is not applied.  Similarly, $2'$ and then $3'$ can each be removed without performing jeu de taquin.  Thus
\begin{equation*}
\lambda^B(w') = (2n-1,2n-3, \ldots, 7,5).
\end{equation*}

Having determined $\lambda^B(w')$, it remains to compute the ratio $f^{\lambda^B(w')}/f^{\lambda^B(\w)}$ via Proposition~\ref{Bhook}.  As in the proof of Theorem~\ref{E(X)thm}, many of the hook-lengths cancel.  Figure~\ref{01shapes} depicts the only boxes in the two shapes where the hook-lengths differ.

\begin{figure}[htbp]
\begin{center}
$\begin{array}{c@{\hspace{.2in}}c}
\epsfig{file=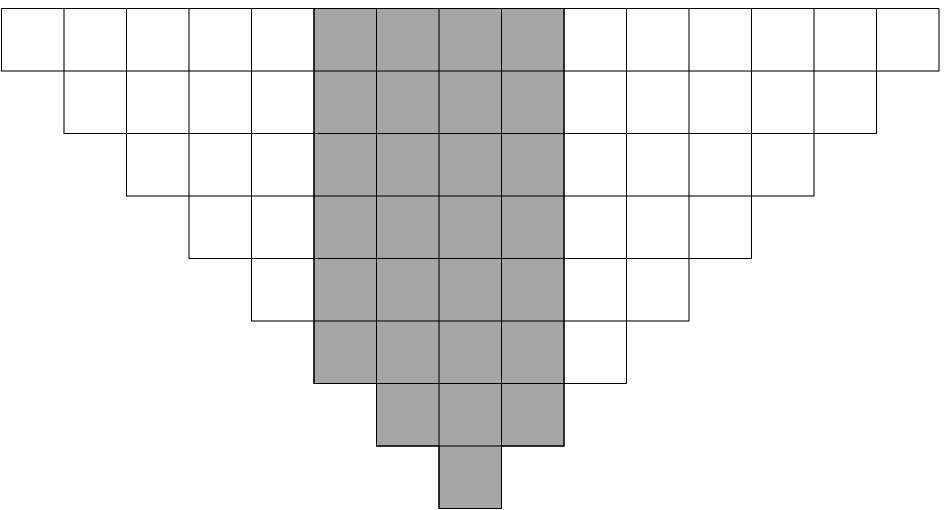, width=2in} & \hspace{.2in} \epsfig{file=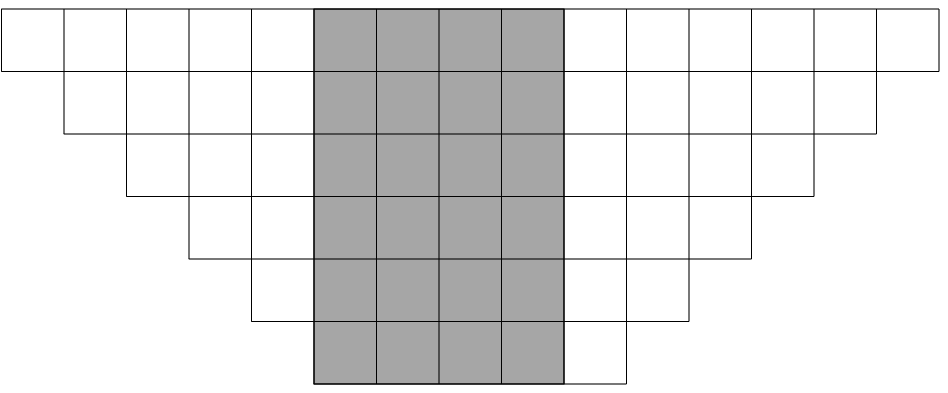, width=2in}
\end{array}$
\end{center}
\caption{The shifted shapes $\lambda^B(\w)$ and $\lambda^B(w')$ for $n = 8$.  The shaded boxes indicate unequal hook-lengths.}\label{01shapes}
\end{figure}

From here it is not hard to compute that
\begin{eqnarray*}
E(\Y) &=& 2(n^2 - 3) \frac{(n^2-4)!}{(n^2)!} \cdot 3 \cdot (2n) \cdot \frac{(2n-2)(2n)(2n+2)}{2\cdot 4\cdot 6}\\
&=& \binom{n^2}{3}^{-1} \cdot (2n) \cdot \frac{(n-1)(n)(n+1)}{1\cdot 2\cdot 3}\\
&=& \frac{2}{n^2-2}.
\end{eqnarray*}

\end{proof}

\section{Acknowledgments}

I am deeply grateful to Vic Reiner for numerous conversations and feedback, and for suggesting the generalization to type $B$.  Thanks also to Dennis Stanton for his comments regarding Dougall's theorem.

\end{document}